\documentclass[11pt]{article}
\usepackage[latin1]{inputenc}
\usepackage[T1]{fontenc}
\usepackage{amsmath,amssymb,amstext}
\usepackage{theorem}
\usepackage{a4wide}
\usepackage{multicol}
\usepackage[english]{babel}
\usepackage{enumerate}
\usepackage{eufrak}

\def\HH{\mathcal{H}}
\def\R{{ \mathbb{R}}}
\newtheorem{theorem}{Theorem}

\newtheorem{lemma}[theorem]{Lemma}

\def\cal{\mathcal}

\renewcommand{\geq}{\geqslant}
\def\leq{\leqslant}

\def\HH{{\cal H}}

\def\cal{\mathcal}
\def\1{{\mathbf{1}}}

\def\Un{{\bf 1}}

\def\sk{{\mathbb{D}}}

\def\oT{{[0{,}T]}}

\def\ffi{{\varphi}}

\def\1{{\mathbf{1}}}
\def\0.5{{\frac{1}{2}}}

\newcommand{\qed}{\nopagebreak\hspace*{\fill}
{\vrule width6pt height6ptdepth0pt}\par}
\renewcommand{\thefootnote}{\fnsymbol{footnote}}

\usepackage{mathrsfs}
\begin{document}

\renewcommand{\thefootnote}{\fnsymbol{footnote}}

\begin{center}
{\large{\bf Parameter estimation for $\alpha$-fractional bridges}}\\~\\
Khalifa Es-Sebaiy\footnote{Institut de Math{\'e}matiques de Bourgogne, Universit{\'e} de
Bourgogne, Dijon,
France.
Email: {\tt khalifasbai@gmail.com}
}
and
Ivan Nourdin\footnote{Institut \'Elie Cartan, Universit\'e Henri Poincar\'e, BP 239, 54506 Vandoeuvre-l\`es-Nancy, France. Email: {\tt inourdin@gmail.com}}\footnote{Supported in part by the (french) ANR grant `Exploration des Chemins Rugueux'}\\
{\it Université de Bourgogne} and {\it Universit\'e Nancy 1}\\~\\
\end{center}
{\small \noindent {\bf Abstract:} 
Let $\alpha,T>0$. We study the asymptotic properties of a least squares estimator
for the parameter $\alpha$ of a fractional bridge defined as
$dX_t=-\alpha\,\frac{X_t}{T-t}\,dt+dB_t$, $0\leq t<T$,
where $B$ is a fractional Brownian motion of Hurst parameter $H>\frac12$.
Depending on the value of $\alpha$, we prove that
we may have strong consistency or not as $t\to T$.
When we have consistency, we obtain
the rate of this convergence as well. Also, 
we compare our results to the (known) case where $B$ is replaced by a standard Brownian motion $W$.\\

{\it It is great pleasure for us to dedicate this paper to our friend David Nualart, in celebration of his 60th birthday and
with all our admiration.}

\section{Introduction}\label{intro}

\
\
\
\
Let $W$ be a standard Brownian motion and let $\alpha$ be a non-negative real parameter.
In recent years, the study of various problems related to the (so-called) $\alpha$-Wiener bridge, 
that is, to the solution $X$ 
to
\begin{equation}\label{eqintro}
X_0=0;\quad dX_t=-\alpha\,\frac{X_t}{T-t}\,dt+dW_t,\quad 0\leq t<T,
\end{equation}
has attracted interest. For a motivation and further references, we refer the reader to Barczy and Pap 
\cite{Barczy,Barczy2}, as well as Mansuy \cite{Mansuy}.
Because (\ref{eqintro}) is linear, it is immediate to solve it explicitely; one then gets the following formula:
\[
X_t = (T-t)^\alpha \int_0^t (T-s)^{-\alpha}dW_s, \quad t\in[0,T),
\]
the integral with respect to $W$ being a Wiener integral.

An example of interesting problem related to $X$ 
is the statistical estimation of $\alpha$ when one observes the whole trajectory of $X$.
A natural candidate  is the maximum likelihood estimator (MLE), which can be easily
computed for this model, due to the specific form of (\ref{eqintro}): one gets
\begin{equation}\label{ac}
\hat{\alpha}_t=
-\left(\int_0^t\frac{X_u}{T-u}\,dX_u\right)\bigg/\left(\int_0^t\frac{X^2_u}{(T-u)^2}\,du\right)
,\quad  t<T.
\end{equation}
In (\ref{ac}), the integral with respect to $X$ must of course be understood in the It\^o sense.
On the other hand, at this stage it is worth noticing that $\hat{\alpha}_t$ coincides with a 
least squares estimator (LSE) as well; indeed,
$\widehat{\alpha}_t$ (formally) minimizes 
\[
\alpha\mapsto \int_0^t \left|\dot{X}_u+\alpha\frac{X_u}{T-u}\right|^2du.
\]
Also, it is worth bearing in mind an alternative formula for $\widehat{\alpha}_t$, 
which is more easily amenable to analysis and which is immediately shown thanks to (\ref{eqintro}):
\begin{eqnarray}\label{ac2}
\alpha - \widehat{\alpha}_t &= &
\left(\int_0^t\frac{X_u}{T-u}\,dW_u\right)\bigg/\left(\int_0^t\frac{X^2_u}{(T-u)^2}\,du\right).
\end{eqnarray}
When dealing with (\ref{ac2}) by means of a semimartingale approach, it is not very difficult to 
check that $\widehat{\alpha}_t$ is indeed a strongly consistent estimator of $\alpha$.
The next step generally consists in studying the second-order approximation. 
Let us describe what is known about this problem: as $t\to T$,
\begin{itemize}
\item 
if $0<\alpha<\frac12$ then \begin{equation}\label{9}
(T-t)^{\alpha-\frac12}\big(\alpha-\widehat{\alpha}_t\big)\overset{\rm law}{\longrightarrow} 
T^{\alpha-\frac12}(1-2\alpha)\times\mathcal{C}(1),\end{equation}
with $\mathcal{C}(1)$ the standard Cauchy distribution,
see \cite[Theorem 2.8]{Barczy3};
\item
if $\alpha=\frac12$ then 
\begin{equation}\label{10}
|\log(T-t)|\big(\alpha-\widehat{\alpha}_t\big)\overset{\rm law}{\longrightarrow}
\frac{\int_0^T W_sdW_s}{\int_0^T W_s^2ds},\end{equation}
see \cite[Theorem 2.5]{Barczy3};
\item
if $\alpha>\frac12$ then \begin{equation}\label{11}
\sqrt{|\log(T-t)|}\big(\alpha-\widehat{\alpha}_t\big)\overset{\rm law}{\longrightarrow} \mathcal{N}(0,2\alpha-1),\end{equation}
see \cite[Theorem 2.11]{Barczy3}.
\end{itemize}
Thus, we have the full picture for the asymptotic behavior of the MLE/LSE associated to $\alpha$-Wiener bridges.\\

In the present paper, our goal is to investigate what happens when, in (\ref{eqintro}), 
the standard Brownian motion $W$ is replaced by a fractional Brownian motion $B$.
More precisely, suppose from now on that
$X=\{X_t\}_{t\in[0,T)}$ 
is the solution
to
\begin{eqnarray}\label{fractionalBridge}
X_0=0;\quad dX_t=-\alpha\,\frac{X_t}{T-t}\,dt+dB_t,\quad 0\leq t<T,
\end{eqnarray} 
where $B$ is a fractional Brownian
motion with known parameter $H$, 
whereas $\alpha>0$ is considered as an unknown parameter.
Although $X$ could have 
been defined for all $H$ in $(0,1)$, for technical reasons and in order to keep the length of our paper within bounds
we restrict ourself
to the case
$H\in(\frac12,1)$ in the sequel.

In order to estimate the unknown parameter $\alpha$ when the whole trajectory of $X$ is
observed, we continue to consider the estimator $\widehat{\alpha}_t$ given by (\ref{ac}).
(It is no longer the MLE, but it is still a LSE.)
Nevertheless, there is a major difference with respect to the standard Brownian motion case.
Indeed, the process $X$ being no longer a semimartingale, in (\ref{ac}) one cannot utilize the It\^o integral 
to integrate
with respect to it.
However, because  $X$ has\footnote{More precisely, we assume throughout the paper that
we work with a suitable $\gamma$-H\"older continuous version of $X$, which is easily shown to 
exist by the
Kolmogorov-Centsov theorem.} $\gamma$-H\"older continuous
paths on $[0,t]$ for all $\gamma\in(\frac12,H)$ and all $t\in[0,T)$,
one can choose, instead, the Young integral (see Section \ref{young} for the main properties of this integral,
notably its chain rule (\ref{chain}) and how (\ref{link}) relies it Skorohod integral).\\
 
Let us now describe the results we prove in the present paper. First,
in Theorem \ref{prop7} we show that
the (strong) consistency of $\widehat{\alpha}_t$ as $t\to T$ holds true if and only if
$\alpha\leq \frac12$.
Then, depending on the precise value of $\alpha\in(0,\frac12]$, we derive the asymptotic behavior of 
the error $\widehat{\alpha}_t-\alpha$. It turns out that, once adequately renormalized,
this error converges either in law or almost surely, to a limit that we are able to compute explicitely.
More specifically,  we show in Theorem \ref{prop8} the following convergences
(below and throughout the paper, $\mathcal{C}(1)$ always stands for the standard Cauchy distribution and 
$\beta(a,b)=\int_0^1 x^{a-1}(1-x)^{b-1}dx$ for the usual
Beta function): as $t\to T$,
\begin{itemize}
\item if $0<\alpha<1-H$ then
\begin{equation}\label{1}
(T-t)^{\alpha-H}\big( \alpha - \widehat{\alpha}_t \big)
\overset{{\rm law}}{\longrightarrow}
T^{\alpha-H}(1-2\alpha)\sqrt{\frac{(H-\alpha)\beta(2-2H-\alpha,2H-1)}{(1-H-\alpha)\beta(1-\alpha,2H-1)}}
\times\mathcal{C}(1);
\end{equation}
\item if $\alpha=1-H$ then 
\begin{equation}\label{2}
\frac{(T-t)^{1- 2H}}{\sqrt{|\log(T-t)|}}
\big( \alpha - \widehat{\alpha}_t \big)
\overset{{\rm law}}{\longrightarrow}
T^{1-2H}
(2H-1)^{\frac32}\sqrt{\frac{2\,\beta(1-H,2H-1)}{\beta(H,2H-1)}}
\times\mathcal{C}(1);
\end{equation}
\item if $1-H<\alpha<\frac12$ then
\begin{equation}\label{3}
(T-t)^{2\alpha-1}\big( \alpha - \widehat{\alpha}_t \big)
\overset{{\rm a.s.}}{\longrightarrow}
(1-2\alpha)\,\int_0^T \frac{dB_u}{(T-u)^{1-\alpha}}\int_0^u \frac{dB_s}{(T-s)^\alpha}\bigg/
\left(\int_0^T \frac{dB_s}{(T-s)^{\alpha}}\right)^2;
\end{equation}
\item if $\alpha= \frac12$ then \begin{equation}\label{4}
|\log(T-t)|\big( \alpha - \widehat{\alpha}_t \big)
\overset{{\rm a.s.}}{\longrightarrow}
\frac12.
\end{equation}
\end{itemize}

When comparing the convergences (\ref{1}) to (\ref{4}) with those arising in the standard Brownian motion case
(that is, (\ref{9}) to (\ref{11})),
we observe a new and interesting phenomenom when the parameter $\alpha$ ranges from $1-H$ to $\frac12$
(of course, this case is immaterial in the standard Brownian motion case).

We hope our proofs of (\ref{1}) to (\ref{4}) to be elementary. Indeed, except maybe  
the link (\ref{link}) between Young and Skorohod integrals, they only involve soft arguments, often
based on the mere derivation of suitable equivalent for some integrals.
In particular, unlike the classical approach (as used, e.g., in \cite{Barczy3}) we stress that, here,  we 
use no tool coming from the semimartingale realm.\\

Before to conclude this introduction, we would like to mention the recent paper \cite{HN} by Hu and Nualart,
which has been a valuable source of inspiration.
More specifically, 
the authors of \cite{HN} study the estimation of the parameter $\alpha>0$ arising in the
fractional Ornstein-Uhlenbeck model, 
defined as $dX_t=-\alpha X_t dt+dB_t$, $t\geq 0$, where $B$ is a fractional Brownian motion
of (known) index $H\in(\frac12,\frac34)$.
They show the strong consistency of a least squares estimator $\widehat{\alpha}_t$ as $t\to\infty$
(with, however, a major difference with respect to us:
they are forced to use Skorohod integral rather than Young integral to define $\widehat{\alpha}_t$, otherwise
$\widehat{\alpha}_t\not\to \alpha$ as $t\to\infty$; unfortunately,
this leads to an impossible-to-simulate estimator, and this is
why they introduce an alternative estimator for $\alpha$.)
They then 
derive the associated rate
of convergence as well, by exhibiting a central limit theorem. Their
calculations are 
of 
completely different nature than ours because, to achieve their goal, the authors of \cite{HN} make use of the fourth moment theorem
of Nualart and Peccati \cite{nualart-peccati}.\\

The rest of our paper is organized as follows.
In Section 2 we introduce the needed material for our study, whereas
Section 3 contains the precise statements and proofs of our results.

\section{Basic notions for fractional Brownian motion}

\
\
\
\
In this section, we briefly recall some basic facts concerning stochastic calculus with
respect to fractional Brownian motion; we refer to
\cite{nualart-book} for further details. Let $B=\{B_t\}_{t\in[0,T]}$ be a
fractional Brownian motion with Hurst parameter $H\in (0,1)$, defined on some probability space $(\Omega,\mathcal{F},P)$.
(Here, and throughout the text, we do assume that $\mathcal{F}$ is the sigma-field generated by $B$.) 
This means
that $B$ is a centered
Gaussian process with the covariance function $E[B_sB_t]=R_H(s,t)$, where%
\begin{align}
\label{cov} R_H(s,t)=\frac{1}2\left (
t^{2H}+s^{2H}-|t-s|^{2H}\right).
\end{align}
If $H=\frac12$, then $B$ is a Brownian motion. From \eqref{cov}, one
can easily see that $E\big[|B_t-B_s|^2\big]=|t-s|^{2H}$, so $B$ has
$\gamma-$H\"{o}lder continuous paths for any $\gamma\in (0,H)$ thanks to 
the Kolmogorov-Centsov theorem.

\subsection{Space of deterministic integrands}\label{sec2.1}
\
\
\
\
We denote by $\mathcal E$ the set of step $\R-$valued functions on
$\oT$. Let $\mathcal H$ be the Hilbert space defined as the
closure of $\mathcal E$ with respect to the scalar product
$$
\left \langle\1_{[0,t]}, \1_{[0,s]}\right \rangle_{\mathcal
H}=R_H(t,s).
$$
We denote by $|\cdot|_{\cal H}$ the associated norm. The mapping
$\displaystyle \1_{[0,t]} \mapsto B_{t}$ can be extended to an
isometry between $\mathcal H$ and the Gaussian space
associated with $B$. We denote this isometry by 
\begin{equation}\label{wienerfrac}
\ffi\mapsto
B(\ffi)=\int_0^T \ffi(s)dB_s.
\end{equation}

When $H  \in (\frac{1}{2}, 1)$, it follows from \cite{PiTa} that the
elements of ${\cal{H}}$ may not be functions but distributions of
negative order. It will be more convenient to work with a subspace
of $\cal H$ which contains only functions. Such a space is the set
$\vert {\cal{H}}\vert$ of all measurable functions $\ffi$ on $[0,T]$
such that
\begin{equation*}
|\ffi|^2_{|\cal H|}:=H (2H -1)\int _{0}^{T} \int _{0}^{T} \vert \ffi(u)
\vert \vert \ffi(v)\vert \vert u-v\vert ^{2H -2} dudv <\infty.
\end{equation*}
If $\ffi,\psi\in|\cal H|$ then 
\begin{equation}\label{iso}
E\big[B(\ffi)B(\psi)\big]=H (2H -1)\int _{0}^{T} \int _{0}^{T} \ffi(u)
\psi(v) \vert u-v\vert ^{2H -2} dudv.
\end{equation}
We know that $(|\cal H|, \langle \cdot,\cdot\rangle_{|\cal H|})$ is a Banach space,
but that $(|\cal H|, \langle\cdot,\cdot\rangle_{\cal H})$ is not
complete (see, e.g., \cite{PiTa}).
We have the dense inclusions
$
L^{2}([0,T]) \subset L^{\frac{1}{H}}([0,T]) \subset
\vert {\cal{H}}\vert \subset {\cal{H}}.
$

\subsection{Malliavin derivative and Skorohod integral}\label{skorohod}
\
\
\
\
Let $\cal S$ be the set of all smooth cylindrical random variables,
which can be expressed as $F = f(B(\phi_1), \ldots, B(\phi_n))$
where $n\geq 1$, $f : \R^n \rightarrow \R$ is a $\mathcal{C}^\infty$-function such that
$f$ and all its derivatives have at most polynomial growth, and $\phi_i\in\cal H$, $i=1,\ldots,n$. The Malliavin derivative of
$F$ with respect to $B$ is the element of $L^2(\Omega, \HH)$ defined by
\[
D_sF\; =\; \sum_{i =1}^n \frac{\partial f}{\partial x_i}(B(\phi_1), \ldots, B(\phi_n))
\phi_i(s), \quad s \in [0,T].
\] 
In particular $D_s B_t =
\Un_{[0,t]}(s)$. As usual, $\sk^{1,2}$ denotes the closure of the set of
smooth random variables with respect to the norm
$$\| F\|_{1,2}^2 \; = \; E[F^2] +
E\big[|DF|_{\HH}^2\big].$$ The Malliavin
derivative $D$ verifies the chain rule: if
$\varphi:\R^n\rightarrow\R$ is $\mathcal{C}^1_b$ and if
$(F_i)_{i=1,\ldots,n}$ is a sequence of elements of $\sk^{1,2}$,
then $\varphi(F_1,\ldots,F_n)\in\sk^{1,2}$ and we have, for any
$s\in [0,T]$,
$$
D_s\varphi(F_1,\ldots,F_n)=\sum_{i=1}^n
\frac{\partial\varphi}{\partial x_i} (F_1,\ldots, F_n)D_sF_i.
$$
The Skorohod integral $\delta$ is the adjoint of the
derivative operator $D$. If a random variable $u \in
L^{2}(\Omega, \HH)$ belongs to the domain  of the Skorohod integral (denoted by ${\rm dom}\delta$), 
that is, if it verifies
$$
|E\langle DF,u\rangle_{\cal H}|\leq c_u\,\sqrt{E[F^2]}\quad\mbox{for any }F\in{\cal S},
$$
then  $\delta(u)$ is defined by the duality
relationship
$$ E[F \delta(u)]= E\big[\langle DF, u \rangle_{\HH}\big], $$ for every $F \in \sk^{1,2}$.
In the sequel, when $t\in[0,T]$ and $u\in{\rm dom} \delta$, 
we shall sometimes write $\int_0^t u_s \delta B_s$ instead of $\delta(u\Un_{[0,t]})$.
If $h\in \HH$, notice moreover that $\int_0^T h_s\delta B_s=\delta(h)=B(h)$.

For every $q\geq 1$, let $\mathcal{H}_{q}$ be the $q$th Wiener chaos
of $B$,
that is, the closed linear subspace of $L^{2}(\Omega)$ generated by the random variables 
$\{H_{q}\left( B\left(
h\right) \right) ,h\in \mathcal{H},\| h\| _{\mathcal{H}}=1\}$, where $H_{q}$ is
the $q$th Hermite polynomial. The mapping $I_{q}(h^{\otimes
q})=H_{q}\left( B\left( h\right) \right) $ provides a linear isometry
between the symmetric tensor product $\mathcal{H}^{\odot q}$ (equipped with the modified norm
$\|\cdot\|_{\mathcal{H}^{\odot q}}=\frac{1}{\sqrt{q!}}
\|\cdot\|_{\mathcal{H}^{\otimes q}}$) and $\mathcal{H}_{q}$. Specifically, for all $f,g\in\mathcal{H}^{\odot q}$ and $q\geq 1$, one has
\begin{equation}\label{isoint}
E\big[I_q(f)I_q(g)\big]=q!\langle f,g\rangle_{\mathcal{H}^{\otimes q}}.
\end{equation}
On the other hand, it is well-known  that any random variable 
$Z$ belonging to $L^2(\Omega)$ admits the following chaotic
expansion:
\begin{equation}\label{chaos}
Z=E[Z]+\sum_{q=1}^\infty I_q(f_q),
\end{equation}
where the series converges in $L^2(\Omega)$ and the kernels $f_q$, belonging to $\HH^{\odot q}$,
 are
uniquely determined by $Z$.

\subsection{Young integral} \label{young}
\
\
\
\
For any $\gamma\in [0,1]$, we denote by $\mathscr{C}^\gamma([0,T])$ the set of 
$\gamma$-H\"older continuous functions, that is, the set
of functions $f:[0,T]\to\R$ such that
\[
|f|_\gamma := \sup_{0\leq s<t\leq T}\frac{|f(t)-f(s)|}{(t-s)^{\gamma}}<\infty.
\]
(Notice the calligraphic difference between a space $\mathscr{C}$ of H\"older continuous functions, 
and  a space $\mathcal{C}$
of continuously differentiable functions!).
We also set $|f|_\infty=\sup_{t\in[0,T]}|f(t)|$, and we equip $\mathscr{C}^\gamma([0,T])$ with the norm \[
\|f\|_\gamma := |f|_\gamma + |f|_\infty.\]
Let $f\in\mathscr{C}^\gamma([0,T])$, and consider the operator $T_f:\mathcal{C}^1([0,T])
\to\mathcal{C}^0([0,T])$ defined as
\[
T_f(g)(t)=\int_0^t f(u)g'(u)du, \quad t\in[0,T].
\]
It can be shown (see, e.g., \cite[Section 2.2]{rv}) that, for any $\beta\in(1-\gamma,1)$, there exists a
constant $C_{\gamma,\beta,T}>0$ depending only on $\gamma$, $\beta$ and $T$ 
such that, for any $g\in\mathscr{C}^\beta([0,T])$,
\[
\left\|\int_0^\cdot f(u)g'(u)du\right\|_\beta \leq C_{\gamma,\beta,T} \|f\|_\gamma \|g\|_\beta.
\]
We deduce that, for any $\gamma\in (0,1)$, any
$f\in\mathscr{C}^\gamma([0,T])$ and any $\beta\in(1-\gamma,1)$, the linear operator 
$T_f:\mathcal{C}^1([0,T])\subset\mathscr{C}^\beta([0,T])\to \mathscr{C}^\beta([0,T])$, 
defined as $T_f(g)=\int_0^\cdot f(u)g'(u)du$, is continuous with respect
to the norm $\|\cdot\|_\beta$. By density, it extends (in an unique way) 
to an operator defined on $\mathscr{C}^\beta$.
As consequence, if $f\in\mathscr{C}^\gamma([0,T])$, if $g\in\mathscr{C}^\beta([0,T])$ and if 
$\gamma+\beta>1$,
then the (so-called) Young integral $\int_0^\cdot f(u)dg(u)$ is (well) defined as being $T_f(g)$.

The Young integral obeys the following chain rule.
Let $\phi:\R^2\to\R$ be a $\mathcal{C}^2$ function, and let $f,g\in\mathscr{C}^\gamma([0,T])$ 
with $\gamma>\frac12$.
Then 
$\int_0^\cdot \frac{\partial\phi}{\partial f}(f(u),g(u))df(u)$
and
$\int_0^\cdot \frac{\partial\phi}{\partial g}(f(u),g(u))dg(u)$
are 
well-defined as Young integrals. Moreover,
for all $t\in[0,T]$, 
\begin{equation}
\phi(f(t),g(t))=\phi(f(0),g(0))+
\int_0^t \frac{\partial\phi}{\partial f}(f(u),g(u))df(u)
+
\int_0^t \frac{\partial\phi}{\partial g}(f(u),g(u))dg(u).
\label{chain}
\end{equation}

\subsection{Link between Young and Skorohod integrals}\label{section-link}
\
\
\
\
Assume $H>\frac12$, and let $u=(u_t)_{t\in[0,T]}$ be 
a process with paths in $\mathscr{C}^\gamma([0,T])$ for some fixed $\gamma>1-H$.
Then, according to the previous section, the integral $\int_0^T u_sdB_s$ exists pathwise in the Young sense.
Suppose moreover that $u_t$ belongs to $\mathbb{D}^{1,2}$ for all $t\in[0,T]$, and that $u$ satisfies
\[
P\left(\int_0^T \int_0^T |D_su_t||t-s|^{2H-2}dsdt<\infty\right)=1.
\]
Then $u\in{\rm dom}\delta$, and we have (see \cite{alos}), for all $t\in[0,T]$:
\begin{equation}\label{link}
\int_0^t u_sdB_s = \int_0^t u_s\delta B_s+H(2H-1)\int_0^t \int_0^t D_su_x|x-s|^{2H-2}dsdx.
\end{equation}
In particular, notice that 
\begin{equation}\label{star}
\int_0^T \ffi_sdB_s = \int_0^T \ffi_s\delta B_s = B(\ffi)
\end{equation}
when $\ffi$ is non-random.

\section{Statement and proofs of our main results}
\
\
\
\
In all this section, we fix a fractional Brownian motion $B$ of Hurst index $H\in(\frac12,1)$,
as well as a parameter $\alpha>0$.
Let us consider the solution $X$ to (\ref{fractionalBridge}).
It is readily checked that we have the following explicit expression for $X_t$:
\begin{equation}\label{explicit}
X_t = (T-t)^\alpha \int_0^t (T-s)^{-\alpha}dB_s, \quad t\in[0,T),
\end{equation}
where the integral can be understood either in the Young sense, or in the Skorohod sense, see indeed (\ref{star}).

For convenience, and because it will play an important role in the forthcoming computations, 
we introduce the following two processes related to $X$: for $t\in[0,T]$,
\begin{eqnarray}
\xi_t&=&\int_0^t (T-s)^{-\alpha}dB_s;\label{xi}\\
\eta_t&=&\int_0^t dB_u(T-u)^{\alpha-1}\int_0^u dB_s(T-s)^{-\alpha}=\int_0^t (T-u)^{\alpha-1}\xi_udB_u\label{eta}.
\end{eqnarray}
In particular, we observe that
\begin{equation}\label{explicit2}
X_t= (T-t)^\alpha \xi_t\quad\mbox{and}\quad 
\int_0^t \frac{X_u}{T-u}dB_u=\eta_t\quad\mbox{for $t\in[0,T)$.}
\end{equation}
When $\alpha$ is between $0$ and $H$ (resp. $1-H$ and $H$), in Lemma \ref{lm3} (resp. Lemma \ref{lm3bis})  
we shall actually show that the process $\xi$ (resp. $\eta$) is well-defined
on the whole interval $[0,T]$ (notice that we could have had a problem at $t=T$), and that
it admits a continuous modification. 
This is why we may and will assume in the sequel, without loss of generality, that $\xi$ (resp. $\eta$) is 
continuous when $0<\alpha<H$ (resp. $1-H<\alpha<H$).

Recall the definition (\ref{ac}) of $\widehat{\alpha}_t$. 
By using (\ref{fractionalBridge}) and then (\ref{explicit2}),
as well as the definitions (\ref{xi}) and (\ref{eta}),
we arrive to the following formula:
\begin{eqnarray*}
\alpha-\widehat{\alpha}_t=
\frac{\int_0^tX_u(T-u)^{-1}dB_u}{\int_0^tX^2_u(T-u)^{-2}ds}
=\frac{\eta_t}
{\int_0^t(T-u)^{2\alpha-2}\xi_u^2du}.
\end{eqnarray*}
Thus, in order to prove the convergences (\ref{1}) to (\ref{4}) of the introduction (that is, our main result!),
we are left to study
the (joint) asymptotic behaviors of $\eta_t$ and $\int_0^t(T-u)^{2\alpha-2}\xi_u^2du$
as $t\to T$.
The asymptotic behavior of $\int_0^t(T-u)^{2\alpha-2}\xi_u^2du$ is rather easy to derive (see Lemma \ref{lm7}),  because
it looks like a convergence {\sl à la} Ces\`aro when $\alpha\leq \frac12$. 
In contrast, the asymptotic behavior of $\eta_t$ is more difficult to obtain, and will depend on the relative position
of $\alpha$ with respect to $1-H$. It is actually the combination of Lemmas \ref{lm1}, \ref{lm3bis}, \ref{lm2},  \ref{lm4},
\ref{lm6} that will allow to derive it for the full range of values of $\alpha$.\\

We are now in position to prove our two main results, that we restate here as theorems for convenience.

\begin{theorem}\label{prop7}
We have
$\widehat{\alpha}_t 
\overset{{\rm prob.}}{\longrightarrow} 
\alpha\wedge \frac12$ as $t\to T$.
When $\alpha<H$ we have almost sure convergence as well.
\end{theorem}

As a corollary, we find 
that $\widehat{\alpha}_t$ is a strong consistent estimator of $\alpha$ if and only if $\alpha\leq\frac12$.
The next result precises the associated rate of convergence in this case.

\begin{theorem}\label{prop8}
Let $G\sim\mathcal{N}(0,1)$ be independent of $B$,
let $\mathcal{C}(1)$ stand for the standard Cauchy distribution, and let $\beta(a,b)=\int_0^1 x^{a-1}(1-x)^{b-1}dx$ 
denote the usual Beta function.
\begin{enumerate}
\item Assume $\alpha\in(0,1-H)$. Then, as $t\to T$,
\begin{eqnarray*}
(T-t)^{\alpha-H}\big( \alpha - \widehat{\alpha}_t \big)
&\,\,\overset{{\rm law}}{\longrightarrow}\,\, &
(1-2\alpha)\sqrt{H(2H-1)\frac{\beta(2-\alpha-2H,2H-1)}{1-H-\alpha}}
\times\frac{G}{\xi_T}\\
&\,\,\overset{{\rm law}}{=}\,\, &
T^{\alpha-H}(1-2\alpha)\sqrt{\frac{(H-\alpha)\beta(2-2H-\alpha,2H-1)}{(1-H-\alpha)\beta(1-\alpha,2H-1)}}
\times\mathcal{C}(1).
\end{eqnarray*}
\item Assume $\alpha=1-H$. Then, as $t\to T$, 
\begin{eqnarray*}
\frac{(T-t)^{1- 2H}}{\sqrt{|\log(T-t)|}}
\big( \alpha - \widehat{\alpha}_t \big)
&\,\,\overset{{\rm law}}{\longrightarrow}\,\,&
(2H-1)^\frac32\sqrt{2H\,\beta(1-H,2H-1)}\times\frac{G}{\xi_T}
\\
&\,\,\overset{{\rm law}}{=}\,\, &
T^{1-2H}
(2H-1)^{\frac32}\sqrt{\frac{2\,\beta(1-H,2H-1)}{\beta(H,2H-1)}}
\times\mathcal{C}(1).
\end{eqnarray*}
\item Assume $\alpha\in\big(1-H,\frac12\big)$. Then, as $t\to T$,
\[
(T-t)^{2\alpha-1}\big( \alpha - \widehat{\alpha}_t \big)
\,\,\overset{{\rm a.s.}}{\longrightarrow}\,\, \frac{(1-2\alpha)\,\eta_T}{(\xi_T)^2}.
\]
\item Assume $\alpha=\frac12$. Then, as $t\to T$,
\[
|\log(T-t)|\big( \alpha - \widehat{\alpha}_t \big)
\,\,\overset{{\rm a.s.}}{\longrightarrow}\,\, \frac12.
\]
\end{enumerate}
\end{theorem}

The rest of this section is devoted to the proofs 
of Theorems \ref{prop7} and \ref{prop8}. Before to be in position to do so, 
we need to state and prove some auxiliary lemmas.
In what follows we use the same symbol $c$ for all
constants whose precise value is not important for our consideration.

\begin{lemma}\label{lm1}
Let $\alpha,\beta\in(0,1)$ be such that $\alpha+\beta<2H$.
Then, for all $T>0$,
\[
\int_0^T ds \,(T-s)^{-\beta}\int_0^T dr\,(T-r)^{-\alpha}|s-r|^{2H-2}
=
\int_0^T ds \,s^{-\beta}\int_0^T dr\,r^{-\alpha}|s-r|^{2H-2}<\infty.
\] 
\end{lemma}
{\it Proof}.
By homogeneity, we first notice that
\[
\int_0^T ds \,s^{-\beta}\int_0^T dr\,r^{-\alpha}|s-r|^{2H-2} = T^{2H-\alpha-\beta}
\int_0^1 ds \,s^{-\beta}\int_0^1 dr\,r^{-\alpha}|s-r|^{2H-2},
\]
so that it is not a loss of generality to assume in the proof that $T=1$.
If $\alpha+1<2H$ then $\int_0^{1/s} r^{-\alpha}|1-r|^{2H-2}dr\leq c s^{-2H+1+\alpha}$,
implying in turn
\[
\int_0^1 ds \,s^{-\beta}\int_0^1 dr\,r^{-\alpha}|s-r|^{2H-2} =\int_0^1 ds \,s^{2H-\alpha-\beta-1}\int_0^{1/s} dr\,r^{-\alpha}|1-r|^{2H-2}\leq c\int_0^1s^{-\beta}ds<\infty.
\]
If $\alpha+1=2H$, then  $\int_0^{1/s} r^{1-2H}|1-r|^{2H-2}dr\leq c (1+|\log s|)$,
implying in turn 
\begin{eqnarray*}
&&\int_0^1 ds \,s^{-\beta}\int_0^1 dr\,r^{-\alpha}|s-r|^{2H-2} 
=\int_0^1 ds \,s^{-\beta}\int_0^1 dr\,r^{1-2H}|s-r|^{2H-2} \\
&=&\int_0^1 ds \,s^{-\beta}
\int_0^{1/s} dr\,r^{1-2H}|1-r|^{2H-2}
\leq c\int_0^1s^{-\beta}  \big(1+   |\log s|\big)ds<\infty.
\end{eqnarray*}
Finally, if $\alpha+1>2H$, then 
\begin{eqnarray*}
\int_0^1 ds \,s^{-\beta}\int_0^1 dr\,r^{-\alpha}|s-r|^{2H-2} 
&=&\int_0^1 ds \,s^{2H-\alpha-\beta-1}\int_0^{1/s} dr\,r^{-\alpha}|1-r|^{2H-2}\\
&\leq& \int_0^1s^{2H-\alpha-\beta-1}ds\times\int_0^{\infty} r^{-\alpha}|1-r|^{2H-2}dr<\infty.
\end{eqnarray*}
\qed

\begin{lemma}\label{lm3}
Assume $\alpha\in(0,H)$. Recall the definition (\ref{xi}) of $\xi_t$.
Then $\xi_T :=\lim_{t\to T} \xi_t$ exists in $L^2$.
Moreover, for all $\varepsilon\in(0,H-\alpha)$, 
the process $\{\xi_t\}_{t\in[0,T]}$ admits a modification 
with $(H-\alpha-\varepsilon)$-H\"older continuous paths, still denoted $\xi$ in the sequel.
In particular,
$\xi_t
\to 
\xi_T$ almost surely as $t\to T$.
\end{lemma}
{\it Proof}.
Because $\alpha<H$, by Lemma \ref{lm1} we have that
$
\int_0^T ds \,s^{-\alpha}\int_0^T du\,u^{-\alpha}|s-u|^{2H-2}<\infty.
$
For all $s\leq t<T$, we thus have, using (\ref{iso}) to get the first equality, 
\begin{eqnarray}
E\left[\left(\xi_t-\xi_s\right)^2\right]
&=&H(2H-1)\int_s^t du(T-u)^{-\alpha}\int_s^t dv(T-v)^{-\alpha}|v-u|^{2H-2}\notag\\
&=&H(2H-1)\int_{T-t}^{T-s} du\,u^{-\alpha}\int_{T-t}^{T-s} dv\,v^{-\alpha}|v-u|^{2H-2}\notag\\
&=&H(2H-1)\int_{0}^{t-s} du\,(u+T-t)^{-\alpha}\int_{0}^{t-s} dv\,(v+T-t)^{-\alpha}|v-u|^{2H-2}\notag\\
&\leq&H(2H-1)\int_{0}^{t-s} du\,u^{-\alpha}\int_{0}^{t-s} dv\,v^{-\alpha}|v-u|^{2H-2}\notag\\
&=&H(2H-1)(t-s)^{2H-2\alpha}\int_{0}^{1} du\,u^{-\alpha}\int_0^{1} dv\,v^{-\alpha}|v-u|^{2H-2}
=c(t-s)^{2H-2\alpha}.\notag
\end{eqnarray}
By the Cauchy criterion, we deduce that $\xi_T :=\lim_{t\to T} \xi_t$ exists in $L^2$.
Moreover, because the process $\xi$ is centered and Gaussian, the Kolmogorov-Centsov theorem applies as well, 
thus leading to the desired conclusion.
\qed

\begin{lemma}\label{lm3bis}
Assume $\alpha\in(1-H,H)$. 
Recall the definition (\ref{eta}) of $\eta_t$.
Then $\eta_T :=\lim_{t\to T} \eta_t$ exists in $L^2$.
Moreover,
there exists $\gamma>0$ such that 
$\{\eta_t\}_{t\in[0,T]}$ admits a modification 
with $\gamma$-H\"older continuous paths, still denoted $\eta$ in the sequel.
In particular,
$\eta_t
\to 
\eta_T$ almost surely as $t\to T$.
\end{lemma}
{\it Proof}.
As a first step, fix $\beta_1,\beta_2\in(1-H,H)$ and let us 
show that there exists $\varepsilon=
\varepsilon(\beta_1,\beta_2,H)>0$ and $c=c(\beta_1,\beta_2,H)>0$ such that, for all $0\leq s\leq t\leq T$, 
\begin{equation}\label{samsam}
\int_{[0,t]\times [s,t]}(T-u)^{-\beta_1}(T-v)^{-\beta_2} |u-v|^{2H-2}dudv\leq c(t-s)^\varepsilon.
\end{equation}
Indeed, we have
\begin{eqnarray*}
&&\int_{[0,t]\times [s,t]}(T-u)^{-\beta_1}(T-v)^{-\beta_2} |u-v|^{2H-2}dudv
=\int_{T-t}^{T}du\,u^{-\beta_1}\int_{T-t}^{T-s}dv\,
v^{-\beta_2} |u-v|^{2H-2}\\
&=&\int_{0}^{t}du(u+T-t)^{-\beta_1}\int_{0}^{t-s}dv(v+T-t)^{-\beta_2} |u-v|^{2H-2}
\leq\int_{0}^{t}du\,u^{-\beta_1}\int_{0}^{t-s}dv\,v^{-\beta_2} |u-v|^{2H-2}\\
&=&\int_{0}^{t-s}du\,u^{-\beta_1}\int_{0}^{t-s}dv\,v^{-\beta_2} |u-v|^{2H-2}
+\int_{t-s}^{t}du\,u^{-\beta_1}\int_{0}^{t-s}dv\,v^{-\beta_2} (u-v)^{2H-2}\\
&=&(t-s)^{2H-\beta_1-\beta_2}\int_{0}^{1}\!\!du\,u^{-\beta_1}\int_{0}^{1}\!\!dv\,v^{-\beta_2} |u-v|^{2H-2}
+\int_{t-s}^{t}\!\!du\,u^{-\beta_1-\beta_2+2H-1}\int_{0}^{(t-s)/u}\!\!\!\!\!dv\,v^{-\beta_2} (1-v)^{2H-2}\\
&\leq&c(t-s)^{2H-\beta_1-\beta_2}+c(t-s)^{1-\beta_2}\int_{t-s}^{t}\!\!du\,u^{-\beta_1}(u-t+s)^{2H-2}\quad\mbox{(see Lemma \ref{lm1} for the first integral
and}\\
&&\hskip9.3cm\mbox{use $1-v\geq 1-\frac{t-s}{u}$ for the second one)}\\
&\leq&c(t-s)^{2H-\beta_1-\beta_2}\left(1+\int_{0}^{s/(t-s)}(w+1)^{-\beta_1}w^{2H-2}dw\right)\\
&=&c(t-s)^{2H-\beta_1-\beta_2}\times\left\{
\begin{array}{ll}
1&\quad\mbox{if $\beta_1>2H-1$}\\
\\
1+|\log(t-s)|&\quad\mbox{if $\beta_1=2H-1$}\\
\\
(t-s)^{-2H+1+\beta_1}&\quad\mbox{if $\beta_1<2H-1$}
\end{array}
\right.
\notag\\
&\leq&
c(t-s)^\varepsilon\quad\mbox{for some $\varepsilon\in (0,1\wedge(2H-\beta_1)-\beta_2)$,}
\end{eqnarray*}
hence (\ref{samsam}) is shown.

Now, let $t<T$. Using (\ref{link}), we can write
\begin{equation}\label{qsd}
\eta_t=\int_0^t \xi_u (T-u)^{\alpha-1}\delta B_u +H(2H-1)\int_0^t du (T-u)^{\alpha-1}\int_0^udv(T-v)^{-\alpha}(u-v)^{2H-2}.
\end{equation}
$\bigg($To have the right to write (\ref{qsd}), according to Section \ref{section-link} we must check that: 
$(i)$ $u\to (T-u)^{\alpha-1}\xi_u$ belongs almost surely to $\mathscr{C}^\gamma([0,t])$ for some $\gamma>1-H$; $(ii)$
$\xi_u\in\mathbb{D}^{1,2}$ for all $u\in[0,t]$, and $(iii)$ $\int_{[0,t]^2}(T-u)^{\alpha-1}|D_v\xi_u|\,|u-v|^{2H-2}dudv<\infty$ almost surely.
To keep the length of this paper within bounds, we will do it completely here, and this will serve as a basis for the proof of
the other instances where a similar verification should have been made as well. The main reason why $(i)$ to $(iii)$
are easy to check is because we are integrating on the compact interval $[0,t]$ with $t$ {\it strictly less than $T$}.

{\it Proof of $(i)$}. Firstly, $u\to (T-u)^{\alpha-1}$ is $\mathcal{C}^\infty$ and bounded on $[0,t]$. Secondly, for $u,v\in[0,t]$
with, say, $u<v$, we have
\begin{eqnarray*}
E[(\xi_u-\xi_v)^2]&=&H(2H-1)\int_u^v dx(T-x)^{-\alpha}\int_u^v dy (T-y)^{-\alpha}|y-x|^{2H-2}\\
&\leq&(T-t)^{-2\alpha} H(2H-1)\int_u^v dx\int_u^v dy|y-x|^{2H-2}
=(T-t)^{-2\alpha}|v-u|^{2H}.
\end{eqnarray*} 
Hence, by combining the Kolmogorov-Centsov theorem with the fact that $\xi$ is Gaussian, we get that (almost) all
the sample paths of $\xi$ are $\theta$-H\"olderian on $[0,t]$ for any $\theta\in(0,H)$. Consequently, by choosing
$\gamma\in(1-H,H)$ (which is possible since $H>1/2$), the proof of $(i)$ is concluded.

{\it Proof of $(ii)$}. This is evident, using the representation (\ref{xi}) of $\xi$ as well as the fact that
$s\to (T-s)^{-\alpha}{\bf 1}_{[0,t]}(s)\in |\HH|$, see Section \ref{sec2.1}.

{\it Proof of $(iii)$}. Here again, it is easy: indeed, we have $D_v\xi_u=(T-v)^{-\alpha}{\bf 1}_{[0,u]}(v)$, so
\[
\int_{[0,t]^2}(T-u)^{\alpha-1}|D_v\xi_u|\,|u-v|^{2H-2}dudv
=\int_{[0,t]^2}(T-u)^{\alpha-1}(T-v)^{-\alpha}\,|u-v|^{2H-2}dudv<\infty.\bigg)
\]

Let us go back to the proof. We deduce from (\ref{qsd}), after setting 
\begin{eqnarray*}
\ffi_t(u,v)=\frac12(T-u\vee v)^{\alpha-1}(T-u\wedge v)^{-\alpha}\,{\bf 1}_{[0,t]^2}(u,v),
\end{eqnarray*}
that
\[
\eta_t=I_2(\ffi_t) +H(2H-1)\int_0^t du (T-u)^{\alpha-1}\int_0^udv(T-v)^{-\alpha}(u-v)^{2H-2}.
\]
Hence, because of (\ref{isoint}),
\begin{eqnarray}
E\left[\left(\eta_t-\eta_s\right)^2\right]
&=&2\|\ffi_t-\ffi_s\|^2_{{\mathcal{H}}^{\otimes 2}}
+H^2(2H-1)^2\left(
\int_s^t du(T-u)^{\alpha-1}\int_0^u dv(T-v)^{-\alpha}(u-v)^{2H-2}
\right)^2.\notag\\
\label{dec1}
\end{eqnarray}
We have, by observing that $\ffi_t-\ffi_s\in|\mathcal{H}|^{\odot 2}$, 
\begin{eqnarray}
&&\|\ffi_t-\ffi_s\|^2_{{\mathcal{H}}^{\otimes 2}}
\notag\\
&=&
H^2(2H-1)^2\int_{[0,T]^4}\big[\ffi_t(u,v)-\ffi_s(u,v)\big]
\big[\ffi_t(x,y)-\ffi_s(x,y)\big]
|u-x|^{2H-2}||v-y|^{2H-2}dudvdxdy\notag\\
&=&
\frac14H^2(2H-1)^2\int_{([0,t]^2\setminus[0,s]^2)^2}(T-u\vee v)^{\alpha-1}(T-x\vee y)^{\alpha-1}
(T-u\wedge v)^{-\alpha}(T-x\wedge y)^{-\alpha}
\notag\\
&&\hskip6cm\times |u-x|^{2H-2}||v-y|^{2H-2}dudvdxdy.\notag
\end{eqnarray}
Taking into account the form of the domain in the previous integral and using that $\ffi_t-\ffi_s$ is symmetric, 
we easily show that
$\|\ffi_t-\ffi_s\|^2_{{\mathcal{H}}^{\otimes 2}}$ is upper bounded (up to constant, and without seeking for sharpness) by a sum of integrals
of the type
\[
\int_{[0,t]\times[s,t]\times[0,T]^2}(T-u)^{-\beta_1}(T-v)^{-\beta_2} 
 (T-x)^{-\beta_3}(T-y)^{-\beta_4}|u-x|^{2H-2}|v-y|^{2H-2}dudvdxdy,
\]
with $\beta_1,\beta_2,\beta_3,\beta_4\in\{\alpha,1-\alpha\}$.
Hence, combining Lemma \ref{lm1} with (\ref{samsam}), we deduce that 
there exists $\varepsilon>0$ small enough and $c>0$ such that, for all $s,t\in[0,T]$,
\begin{equation}\label{sum1}
\|\ffi_t-\ffi_s\|^2_{{\mathcal{H}}^{\otimes 2}}\leq c|t-s|^\varepsilon.
\end{equation}

On the other hand, we can write, for all $s\leq t<T$,
\begin{eqnarray}
&&\int_s^t du(T-u)^{\alpha-1}\int_0^u dv(T-v)^{-\alpha}(u-v)^{2H-2}\notag\\
&=&\int_{T-t}^{T-s} du\,u^{\alpha-1}\int_u^T dv\,v^{-\alpha}(v-u)^{2H-2}\notag\\
&=&\int_{0}^{t-s} du(u+T-t)^{\alpha-1}\int_u^t dv(v+T-t)^{-\alpha}(v-u)^{2H-2}\notag\\
&\leq &\int_{0}^{t-s} du\,u^{\alpha-1}\int_u^T dv\,v^{-\alpha}(v-u)^{2H-2}\notag\\
&=&(t-s)^{2H-1}\int_{0}^{1} du\,u^{\alpha-1}\int_u^{\frac{T}{t-s}} dv\,v^{-\alpha}(v-u)^{2H-2}\notag\\
&= &(t-s)^{2H-1}\int_{0}^{1} du\,u^{2H-2} \int_1^{\frac{T}{(t-s)u}} dv\,v^{-\alpha}(v-1)^{2H-2}.\label{grf}
\end{eqnarray}
Let us consider three cases.
Assume first that $\alpha>2H-1$: in this case, 
\[
\int_1^{\frac{T}{(t-s)u}} v^{-\alpha}(v-1)^{2H-2}dv\leq \int_1^{\infty} v^{-\alpha}(v-1)^{2H-2}dv<\infty;
\]
leading, thanks to (\ref{grf}), to
\[
\int_s^t du(T-u)^{\alpha-1}\int_0^u dv(T-v)^{-\alpha}(u-v)^{2H-2}\leq c(t-s)^{2H-1}.
\]
The second case is when $\alpha=2H-1$: we then have 
\[
\int_1^{\frac{T}{(t-s)u}} v^{-\alpha}(v-1)^{2H-2}dv
\leq c\big(1+|\log(t-s)|+|\log u|\big)
\]
so that, by (\ref{grf}),
\[
\int_s^t du(T-u)^{\alpha-1}\int_0^u dv(T-v)^{-\alpha}(u-v)^{2H-2}\leq c(t-s)^{2H-1}\big(1+|\log(t-s)|\big).
\]
Finally, the third case is when $\alpha<2H-1$: in this case, 
\[
\int_1^{\frac{T}{(t-s)u}} v^{-\alpha}(v-1)^{2H-2}dv\leq c (t-s)^{\alpha-2H+1} u^{\alpha-2H+1};
\]
so that, by (\ref{grf}),
\[
\int_s^t du(T-u)^{\alpha-1}\int_0^u dv(T-v)^{-\alpha}(u-v)^{2H-2}\leq c(t-s)^{\alpha}.
\]
To summarize, we have shown that there exists $c>0$ such that, for all $s,t\in[0,T]$,
\begin{equation}\label{sum2}
\int_s^t du(T-u)^{\alpha-1}\int_0^u dv(T-v)^{-\alpha}(u-v)^{2H-2}\leq
c\big(1+|\log(|t-s|)|{\bf 1}_{\{\alpha=2H-1\}}\big)|t-s|^{(2H-1)\wedge \alpha}.
\end{equation}

By inserting (\ref{sum1}) and (\ref{sum2}) into (\ref{dec1}), 
we finally get that there exists $\varepsilon>0$ small enough and $c>0$ such that, for all $s,t\in[0,T]$,
\[
E\left[\left(\eta_t-\eta_s\right)^2\right]\leq c|t-s|^{\varepsilon}.
\]
By the Cauchy criterion, we deduce that $\eta_T :=\lim_{t\to T} \eta_t$ exists in $L^2$.
Moreover, because
$\eta_t-\eta_s-E[\eta_t]+E[\eta_s]$
belongs to the second Wiener chaos of $B$ (where all the $L^p$ norms are equivalent),  
the Kolmogorov-Centsov theorem applies as well, 
thus leading to the desired conclusion.
\qed

\begin{lemma}\label{lm2}
Recall the definition (\ref{eta}) of $\eta_t$.
For any $t\in[0,T)$, we have 
\begin{eqnarray*}
\eta_t
&=& 
\int_0^t (T-u)^{\alpha-1}dB_u\times \int_0^t (T-s)^{-\alpha}dB_s  
-
\int_0^t \delta B_s\,(T-s)^{-\alpha}\int_0^s \delta B_u\,(T-u)^{\alpha-1}\\
&&-H(2H-1)\int_0^t ds\,(T-s)^{-\alpha}\int_0^s du\,(T-u)^{\alpha-1}(s-u)^{2H-2}.
\end{eqnarray*}
\end{lemma}
{\it Proof}. 
Fix $t\in[0,T)$. 
Applying the change of variable formula (\ref{chain}) to the right-hand side of the first equality in 
(\ref{eta}) leads to
\begin{equation}\label{qsd2}
\eta_t =
\int_0^t (T-u)^{\alpha-1}dB_u\times \int_0^t (T-s)^{-\alpha}dB_s  
-
\int_0^t dB_s\,(T-s)^{-\alpha}\int_0^s dB_u\,(T-u)^{\alpha-1}.
\end{equation}
On the other hand, by (\ref{link}) we have that
\begin{eqnarray}
&&\int_0^t dB_s\,(T-s)^{-\alpha}\int_0^s dB_u\,(T-u)^{\alpha-1}\label{qsd3}\\
&\!\!=&\!\!\!\int_0^t \delta B_s\,(T-s)^{-\alpha}\int_0^s \delta B_u\,(T-u)^{\alpha-1}
+H(2H-1)\int_0^t ds(T-s)^{-\alpha}\int_0^s du(T-u)^{\alpha-1}(s-u)^{2H-2}.\notag
\end{eqnarray}
The desired conclusion follows.
(We omit the justification of (\ref{qsd2}) and (\ref{qsd3}) because it suffices to proceed as in the proof (\ref{qsd}).) 
\qed

\begin{lemma}\label{lm4}
Let $\beta(a,b)=\int_0^1 x^{a-1}(1-x)^{b-1}dx$ denote the usual Beta function,
let $Z$ be any $\sigma\{B\}$-measurable random variable satisfying $P(Z<\infty)=1$, and
let $G\sim\mathcal{N}(0,1)$ be independent of $B$.
\begin{enumerate}
\item Assume $\alpha\in(0,1-H)$.
Then, as $t\to T$,
\begin{equation}\label{robot}
\left(
Z,(T-t)^{1-H-\alpha}
\int_0^t (T-u)^{\alpha-1}dB_u\right)
\overset{\rm law}{\longrightarrow} 
\left(Z,\sqrt{H(2H-1)\frac{\beta(2-\alpha-2H,2H-1)}{1-H-\alpha}
}\,G\right).
\end{equation}
\item Assume $\alpha=1-H$.
Then, as $t\to T$,
\begin{equation}\label{robot2}
\left(Z,\frac{1}{\sqrt{|\log(T-t)|}}
\int_0^t (T-u)^{-H}dB_u\right)
\overset{\rm law}{\longrightarrow} 
\left(Z,\sqrt{2H(2H-1)\beta(1-H,2H-1)}\,G\right).
\end{equation}
\end{enumerate}
\end{lemma}
{\it Proof}. By a standard approximation procedure, we first notice that it is not a loss of generality to assume that $Z$
belongs to $L^2(\Omega)$ (using e.g. that  $Z\,{\bf 1}_{\{|Z|\leq n\}}\overset{\rm a.s.}{\longrightarrow} Z$
as $n\to\infty$).\\

1. Set $N=\sqrt{H(2H-1)\frac{\beta(2-\alpha-2H,2H-1)}{1-H-\alpha}
}\,G$.
For any $d\geq 1$ and any $s_1,\ldots,s_d\in[0,T)$, we shall prove that
\begin{equation}\label{toto}
\left(B_{s_1},\ldots,B_{s_d},(T-t)^{1-H-\alpha}
\int_0^t (T-u)^{\alpha-1}dB_u\right)
\overset{{\rm law}}{\longrightarrow} 
\big(B_{s_1},\ldots,B_{s_d},N\big)\quad\mbox{as $t\to T$}.
\end{equation}
Suppose for a moment that $(\ref{toto})$ has been shown, and let us proceed with the proof of (\ref{robot}).
By the very construction of $\HH$ and by reasoning by approximation, we deduce that, for any $l\geq 1$
and any $h_1,\ldots,h_l\in\HH$ with unit norms,  
\[
\left(B(h_1),\ldots,B(h_l),(T-t)^{1-H-\alpha}
\int_0^t (T-u)^{\alpha-1}dB_u\right)
\overset{{\rm law}}{\longrightarrow} 
\big(B(h_1),\ldots,B(h_l),N\big)\quad\mbox{as $t\to T$}.
\]
This implies that, for any $l\geq 1$, any $h_1,\ldots,h_l\in\HH$ with unit norms
and any integers $q_1,\ldots,q_l\geq 0$,
\begin{eqnarray*}
&&\left(H_{q_1}(B(h_1)),\ldots,H_{q_l}(B(h_l)),(T-t)^{1-H-\alpha}
\int_0^t (T-u)^{\alpha-1}dB_u\right)\\
&\overset{{\rm law}}{\longrightarrow}&
\big(H_{q_1}(B(h_1)),\ldots,H_{q_l}(B(h_l)),N\big)\quad\mbox{as $t\to T$},
\end{eqnarray*}
with $H_q$ the $q$th Hermite polynomial.
Using now the very definition of the Wiener chaoses and by reasoning by approximation once again, we deduce that, for any $l\geq 1$,
any integers $q_1,\ldots,q_l\geq 0$ and any $f_1\in\HH^{\odot q_1}, \ldots,f_l\in\HH^{\odot q_l}$,
\begin{eqnarray*}
\left(I_{q_1}(f_1),\ldots,I_{q_l}(f_l),(T-t)^{1-H-\alpha}
\int_0^t (T-u)^{\alpha-1}dB_u\right)
\overset{{\rm law}}{\longrightarrow}
\big(I_{q_1}(f_1),\ldots,I_{q_l}(f_l),N\big)\quad\mbox{as $t\to T$}.
\end{eqnarray*}
Thus, for any random variable $F\in L^2(\Omega)$ with a {\sl finite} chaotic decomposition,
we have
\begin{equation}\label{aqw}
\left(F,(T-t)^{1-H-\alpha}
\int_0^t (T-u)^{\alpha-1}dB_u\right)
\overset{{\rm law}}{\longrightarrow} 
\big(F,N\big)\quad\mbox{as $t\to T$}.
\end{equation}
To conclude, let us consider the chaotic decomposition (\ref{chaos}) of $Z$.
By applying (\ref{aqw}) to $F=E[Z]+\sum_{q=1}^n I_q(f_q)$ and then letting $n\to\infty$, 
we finally deduce that (\ref{robot}) holds true.

Now, let us proceed with the proof of (\ref{toto}).
Because the left-hand side of (\ref{toto}) is a Gaussian vector, to get (\ref{toto}) it is sufficient
to check the convergence of covariance matrices.
Let us first compute the limiting  variance of $(T-t)^{1-H-\alpha}
\int_0^t (T-u)^{\alpha-1}dB_u$ as $t\to T$. By (\ref{iso}), for any $t\in[0,T)$ we have
\begin{eqnarray*}
&&E\left[\left((T-t)^{1-H-\alpha}\int_0^t (T-u)^{\alpha-1}dB_u\right)^2\right]\\
&=&H(2H-1)(T-t)^{2-2H-2\alpha}\int_0^t ds(T-s)^{\alpha-1}\int_0^t du(T-u)^{\alpha-1}|s-u|^{2H-2}\\
&=&H(2H-1)(T-t)^{2-2H-2\alpha}\int_{T-t}^T ds\,s^{\alpha-1}\int_{T-t}^T du\,u^{\alpha-1}|s-u|^{2H-2}\\
&=&H(2H-1)\int_{1}^{\frac{T}{T-t}} ds\,s^{\alpha-1}\int_{1}^{\frac{T}{T-t}} du\,u^{\alpha-1}|s-u|^{2H-2}\\
&\to&H(2H-1)\int_{1}^{\infty} ds\,s^{\alpha-1}\int_{1}^{\infty} du\,u^{\alpha-1}|s-u|^{2H-2}\quad\mbox{as $t\to T$},
\end{eqnarray*}
with
\begin{eqnarray*}
&&\int_{1}^{\infty} ds\,s^{\alpha-1}\int_{1}^{\infty} du\,u^{\alpha-1}|s-u|^{2H-2}
=\int_{1}^{\infty} ds\,s^{2\alpha+2H-3}\int_{1/s}^{\infty} du\,u^{\alpha-1}|1-u|^{2H-2}
\\
&=&\int_{1}^{\infty} s^{2\alpha+2H-3}ds\int_{1}^{\infty} u^{\alpha-1}(u-1)^{2H-2}du
+\int_{1}^{\infty} ds\,s^{2\alpha+2H-3}\int_{1/s}^{1} du\,u^{\alpha-1}(1-u)^{2H-2}\\
&=&\frac{\beta(2-\alpha-2H,2H-1)}{2(1-H-\alpha)} 
+\int_{0}^{1}du \,u^{\alpha-1}(1-u)^{2H-2}\int_{1/u}^\infty ds\,s^{2\alpha+2H-3}\\
&=&\frac{\beta(2-\alpha-2H,2H-1)}{1-H-\alpha}.
\end{eqnarray*}
Thus,
\[
\lim_{t\to T}E\left[\left((T-t)^{1-H-\alpha}\int_0^t (T-u)^{\alpha-1}dB_u\right)^2\right]=
\frac{H(2H-1)}{1-H-\alpha}\,\beta(2-\alpha-2H,2H-1).
\]
On the other hand, by (\ref{iso}) we have, for any $v<t<T$,
\begin{eqnarray*}
&&E\left[B_v\times (T-t)^{1-H-\alpha}\int_0^t (T-u)^{\alpha-1}dB_u\right]\\
&=&H(2H-1)(T-t)^{1-H-\alpha}\int_0^t du\,(T-u)^{\alpha-1}\int_0^v ds\,|u-s|^{2H-2}\\
&=&H(T-t)^{1-H-\alpha}\int_0^t (T-u)^{\alpha-1}\big(u^{2H-1} + {\rm sign}(v-u)\times |v-u|^{2H-1}\big)du\\
&\to&0\quad\mbox{as $t\to T$},
\end{eqnarray*}
because $\int_0^T (T-u)^{\alpha-1}\big(u^{2H-1} + {\rm sign}(v-u)\times |v-u|^{2H-1}\big)du<\infty$.
Convergence (\ref{toto}) is then shown, and (\ref{robot}) follows.\\

2. By (\ref{iso}), for any $t\in[0\vee(T-1),T)$  we have
\begin{eqnarray*}
&&E\left[\left( \frac{1}{\sqrt{|\log(T-t)|}}\int_0^t(T-u)^{-H}dB_u\right)^2\right]\\
&=&\frac{H(2H-1)}{{|\log(T-t)|}}
\int_0^tds(T-s)^{-H}\int_0^tdu(T-u)^{-H}|s-u|^{2H-2}\\
&=&\frac{H(2H-1)}{{|\log(T-t)|}}
\int_{T-t}^Tds\,s^{-H}\int_{T-t}^Tdu\,u^{-H}|s-u|^{2H-2}\\
&=&\frac{2H(2H-1)}{{|\log(T-t)|}}
\int_{T-t}^Tds\,s^{-H}\int_{T-t}^sdu\,u^{-H}(s-u)^{2H-2}\\
&=&\frac{2H(2H-1)}{{|\log(T-t)|}}
\int_{T-t}^T  \frac{ds}{s}  \int_{\frac{T-t}{s}}^1du\,u^{-H}(1-u)^{2H-2}\\
&=&\frac{2H(2H-1)}{{|\log(T-t)|}}
\int_{\frac{T-t}{T}}^1du\,u^{-H}(1-u)^{2H-2}\int_{\frac{T-t}{u}}^T \frac{ds}{s}\\
&=&2H(2H-1)
\int_{\frac{T-t}{T}}^1du\,u^{-H}(1-u)^{2H-2}\left(1+\frac{\log(Tu)}{|\log(T-t)|}\right).
\end{eqnarray*}
Because 
$\int_{0}^1 |\log(Tu)|u^{-H}(1-u)^{2H-2}du<\infty$, we get that 
\[
E\left[\left( \frac{1}{\sqrt{|\log(T-t)|}}\int_0^t(T-s)^{-H}dB_s\right)^2\right]\to 
2H(2H-1)\beta(1-H,2H-1)\quad\mbox{as $t\to T$.}
\] 
On the other hand, fix $v\in[0,T)$. For all $t\in[0\vee(T-1),T)$, using (\ref{iso}) we can write
\begin{eqnarray*}
&&E\left[B_v\times\frac{1}{\sqrt{|\log(T-t)|}}\int_0^t(T-u)^{-H}dB_u\right]
\\
&=&\frac{H(2H-1)}{\sqrt{|\log(T-t)|}}\int_0^tdu\,(T-u)^{-H}\int_0^vds\,|u-s|^{2H-2}\\
&=&\frac{H}{\sqrt{|\log(T-t)|}}
\int_0^T (T-u)^{\alpha-1}\big(u^{2H-1} + {\rm sign}(v-u)\times |v-u|^{2H-1}\big)du\\
&\longrightarrow& 0 \quad \mbox{as $t\to T$,}
\end{eqnarray*}
because $\int_0^T (T-u)^{\alpha-1}\big(u^{2H-1} + {\rm sign}(v-u)\times |v-u|^{2H-1}\big)du<\infty$.
Thus, we have shown that,
for any $d\geq 1$ and any $s_1,\ldots,s_d\in[0,T)$,
\begin{equation}\label{toto2}
\left(B_{s_1},\ldots,B_{s_d},(T-t)^{1-H-\alpha}
\int_0^t (T-u)^{\alpha-1}dB_u\right)
\overset{{\rm law}}{\longrightarrow} 
\bigg(B_{s_1},\ldots,B_{s_d},\sqrt{2H(2H-1)\beta(1-H,2H-1)}\,G\bigg)
\end{equation}
as $t\to T$.
Finally, the same reasoning as in point 1 above allows to go from (\ref{toto2}) to (\ref{robot2}).
The proof of the lemma is concluded.
\qed

\begin{lemma}\label{lm6}
Assume $\alpha\in(0,1-H]$.
Then, as $t\to T$,
\[
\limsup_{t\to T}E\left[\left(
\int_0^t \delta B_u\,(T-u)^{-\alpha}\int_0^s \delta B_v\,(T-v)^{\alpha-1}\right)^2\right]<\infty.
\]
\end{lemma}
{\it Proof}.
Set $\phi_t(u,v)=\frac12
(T-u\vee v)^{-\alpha}(T-u\wedge v)^{\alpha-1}
{\bf 1}_{[0,t]^2}(u,v)$.
We have $\phi_t\in|\HH|^{\odot 2}$ and 
$\int_0^t \delta B_u\,(T-u)^{-\alpha}\int_0^u \delta B_v\,(T-v)^{\alpha-1}=I_2(\phi_t)$, so that
\begin{eqnarray*}
&&\limsup_{t\to T}E\left[\left(
\int_0^t \delta B_u\,(T-u)^{-\alpha}\int_0^u \delta B_v\,(T-v)^{\alpha-1}
\right)^2
\right]=2\,\limsup_{t\to T}\|\phi_t\|^2_{\HH^{\otimes 2}}\\
&=&2H^2(2H-1)^2\,\limsup_{t\to T}\int_{[0,T]^4} \phi_t(u,v)\phi_t(x,y)|u-x|^{2H-2}|v-y|^{2H-2}dudvdxdy\\
&=&\frac12H^2(2H-1)^2\int_{[0,T]^4} (T-u\vee v)^{-\alpha}(T-u\wedge v)^{\alpha-1}
(T-x\vee y)^{-\alpha}(T-x\wedge y)^{\alpha-1}\\
&&\hskip6cm\times|u-x|^{2H-2}|v-y|^{2H-2}dudvdxdy\\
&=&2H^2(2H-1)^2\int_0^T du\,(T-u)^{-\alpha}\int_0^T dx\,(T-x)^{-\alpha}|u-x|^{2H-2} \int_0^u dv\,(T-v)^{\alpha-1} \\
&&\hskip6cm\times \int_0^x dv \, (T-y)^{\alpha-1}|v-y|^{2H-2}\\
&=&2H^2(2H-1)^2\int_{0}^T du\,u^{-\alpha} \int_{0}^T dx\,x^{-\alpha}
 |u-x|^{2H-2} \int_u^T dv\,v^{\alpha-1} \int_x^T dy \, y^{\alpha-1}|v-y|^{2H-2}\\
&=&2H^2(2H-1)^2
\int_{0}^T du\,u^{-\alpha} \int_{0}^T dx\,x^{-\alpha}|u-x|^{2H-2}\int_u^T dv\,v^{2H+2\alpha-3} 
\int_{x/v}^{T/v} dy \, y^{\alpha-1}|1-y|^{2H-2}.
\end{eqnarray*}
Because $\alpha\leq 1-H$ and $H<1$, we have $\alpha<2-2H$, 
so that \[
\int_{x/v}^{T/v} y^{\alpha-1}|1-y|^{2H-2}dy\leq 
\int_{0}^{\infty} y^{\alpha-1}|1-y|^{2H-2}dy <\infty.
\]
Moreover, because $2H+2\alpha-3\leq -1$ due to our assumption on $\alpha$, we have
\[
\int_u^T dv\,v^{2H+2\alpha-3} 
\leq c\left\{\begin{array}{ll}
u^{2H+2\alpha-2}&\quad\mbox{if $\alpha<1-H$}\\
1+|\log u|&\quad\mbox{if $\alpha=1-H$}
\end{array}\right..
\]
Consequently, if $\alpha=1-H$, then
\begin{eqnarray*}
&&\int_{0}^T du\,u^{-\alpha} \int_{0}^T dx\,x^{-\alpha}|u-x|^{2H-2}\int_u^T dv\,v^{2H+2\alpha-3} 
\int_{x/v}^{T/v} dy \, y^{\alpha-1}|1-y|^{2H-2}\\
&\leq&c\int_{0}^T du\,u^{H-1}\big(1+|\log u|\big) \int_{0}^T dx\,x^{H-1}|u-x|^{2H-2}\\
&=&c\int_{0}^T du\,u^{4H-3}\big(1+|\log u|\big) \int_{0}^{T/u} dx\,x^{H-1}|1-x|^{2H-2}\\
&\leq&c\int_{0}^T du \,u^{4H-3}\big(1+|\log u|\big)
\times
\left\{
\begin{array}{ll}
1&\quad\mbox{if $H<\frac23$}\\
\\
1+|\log u|&\quad\mbox{if $H=\frac23$}\\
\\
u^{2-3H}&\quad\mbox{if $H>\frac23$}\\
\end{array}\right.\\
&<&\infty,
\end{eqnarray*}
and the proof is concluded in this case.
Assume now that $\alpha<1-H$. Then
\begin{eqnarray*}
&&\int_{0}^T du\,u^{-\alpha} \int_{0}^T dx\,x^{-\alpha}|u-x|^{2H-2}\int_u^T dv\,v^{2H+2\alpha-3} 
\int_{x/v}^{T/v} dy \, y^{\alpha-1}|1-y|^{2H-2}\\
&\leq&c\int_{0}^T du\,u^{2H+\alpha-2} \int_{0}^T dx\,x^{-\alpha}|u-x|^{2H-2}
=c\int_{0}^T du\,u^{4H-3} \int_{0}^{T/u} dx\,x^{-\alpha}|1-x|^{2H-2}.
\end{eqnarray*}
Let us distinguish three different cases.
First, if $\alpha<2H-1$ then
\[
\int_{0}^T du\,u^{4H-3} \int_{0}^{T/u} dx\,x^{-\alpha}|1-x|^{2H-2}
\leq c\int_{0}^T u^{2H-2+\alpha} du<\infty.
\]
Second, if $\alpha=2H-1$ then 
\begin{eqnarray*}
\int_{0}^T du\,u^{4H-3} \int_{0}^{T/u} dx\,x^{-\alpha}|1-x|^{2H-2}
&=&
\int_{0}^T du\,u^{4H-3} \int_{0}^{T/u} dx\,x^{1-2H}|1-x|^{2H-2}\\
&\leq& c\int_{0}^T u^{4H-3}\big(1+|\log u|\big) du<\infty.
\end{eqnarray*}
Third, if $\alpha>2H-1$ then 
\[
\int_{0}^T du\,u^{4H-3} \int_{0}^{T/u} dx\,x^{-\alpha}|1-x|^{2H-2}\leq \int_{0}^T u^{4H-3} du
\int_{0}^{\infty} x^{-\alpha}|1-x|^{2H-2}dx<\infty.
\]
Thus, in all the possible cases we see that 
$\limsup_{t\to T}E\left[\left(\int_0^t \delta B_u \,(T-u)^{-\alpha}\int_0^u 
\delta B_v\,(T-v)^{\alpha-1}\right)^2\right]$ is finite, and the proof of the lemma is done.\qed

\begin{lemma}\label{lm7}
Assume $\alpha\in(0,H)$, and recall the definition (\ref{xi}) of $\xi_t$.
Then, as $t\to T$:
\begin{enumerate}
\item if $0<\alpha<\frac12$, then \[
(T-t)^{1-2\alpha}\int_0^t \xi_s^2(T-s)^{2\alpha-2}\,ds
\,\,\overset{\rm a.s.}{\to}\,\, \frac{\xi_T^2}{1-2\alpha};\]
\item if $\alpha=\frac12$, then \[\frac{1}{|\log(T-t)|}\int_0^t \frac{\xi_s^2}{T-s}\,ds
\,\,\overset{\rm a.s.}{\to}\,\, \xi_T^2;\]
\item if $\frac12<\alpha<H$, then \[\int_0^t \xi_s^2(T-s)^{2\alpha-2}\,ds
\,\,\overset{\rm a.s.}{\to}\,\,  \int_0^T \xi_s^2(T-s)^{2\alpha-2}\,ds<\infty.\]
\end{enumerate}
\end{lemma}
{\it Proof}. 1. Using the $\big(\frac{H}2-\frac{\alpha}2\big)$-H\"olderianity of
$\xi$ (Lemma \ref{lm3}), we can write
\begin{eqnarray*}
&&\left|(T-t)^{1-2\alpha}\int_0^t \xi_s^2(T-s)^{2\alpha-2}\,ds -\frac{\xi_T^2}{1-2\alpha}\right|\\
&\leq&(T-t)^{1-2\alpha}\int_0^t \big|\xi_s^2-\xi_T^2\big|(T-s)^{2\alpha-2}\,ds 
+(T-t)^{1-2\alpha}\frac{T^{2\alpha-1}}{1-2\alpha}\xi_T^2\\
&\leq&c|\xi|_\infty(T-t)^{1-2\alpha}\int_0^t (T-s)^{\frac{H}2+\frac{3\alpha}2-2}\,ds 
+(T-t)^{1-2\alpha}\frac{T^{2\alpha-1}}{1-2\alpha}\xi_T^2\\
&\leq&c|\xi|_\infty\big( (T-t)^{\frac{H}2-\frac{\alpha}2} + (T-t)^{1-2\alpha}T^{\frac{H}2+\frac{3\alpha}2-1}\big) 
+(T-t)^{1-2\alpha}\frac{T^{2\alpha-1}}{1-2\alpha}\xi_T^2\\
&\to& 0\quad\mbox{almost surely as $t\to T$.}
\end{eqnarray*}
2. 
Using the $(\frac{H}2-\frac14)$-H\"{o}lderianity of $\xi$ (Lemma \ref{lm3}), we can write
\begin{eqnarray*}
&& \left|\frac{1}{|\log(T-t)|}\int_0^t\frac{\xi_s^2}{T-s} ds- \xi_T^2 \right|\\
&\leq&\frac{1}{|\log(T-t)|}\int_0^t\frac{\left|\xi_s^2-\xi_T^2\right|}{T-s} ds+\frac{\log(T)}{|\log(T-t)|} \xi_T^2\\
&\leq&\frac{c|\xi|_{\infty}}{|\log(T-t)|}\int_0^t (T-s)^{\frac{H}{2}-\frac{5}{4}}  ds
+\frac{\log(T)}{|\log(T-t)|} \xi_T^2\\
&\leq&\frac{c|\xi|_{\infty}}{|\log(T-t)|}\big(T^{\frac{H}2-\frac14}+(T-t)^{\frac{H}2-\frac14}\big)+\frac{\log(T)}{|\log(T-t)|} \xi_T^2\\
&\to&0 \quad\mbox{almost surely as $t\to T$.}
\end{eqnarray*}
3. 
By Lemma \ref{lm3}, the process $\xi$ is continuous on $[0,T]$, hence integrable.
Moreover, $s\mapsto (T-s)^{2\alpha-2}$ is integrable at $s=T$ because $\alpha>\frac12$.
The convergence in point 3 is then clear, with a finite limit.

\qed

\medskip

We are now ready to prove Theorems \ref{prop7} and \ref{prop8}.

\bigskip

\noindent
{\it Proof of Theorem \ref{prop7}}.
Fix $\alpha>0$. Thanks to the change of variable formula (\ref{chain})
(which can be well applied here, as is easily shown by proceeding as in the proof (\ref{qsd})), we can write, for any $t\in[0,T)$:
\begin{eqnarray*}
\frac12(T-t)^{2\alpha-1}\xi_t^2&=&
\frac{1-2\alpha}2\int_0^t (T-u)^{2\alpha-2}\xi_u^2du + \int_0^t (T-u)^{2\alpha-1}\xi_u d\xi_u\\
&=&
\frac{1-2\alpha}2\int_0^t (T-u)^{2\alpha-2}\xi_u^2du + \eta_t,
\end{eqnarray*}
so that 
\begin{equation}\label{sothat}
 \alpha - \widehat{\alpha}_t =
\frac{\xi_t^2}{2
(T-t)^{1-2\alpha}\int_0^t \xi_u^2(T-u)^{2\alpha-2}du
}
+\alpha-\frac12.
\end{equation}

When $\alpha\in(0,\frac12)$, we have $(T-t)^{1-2\alpha}\int_0^t \xi_u^2(T-u)^{2\alpha-2}du
\overset{\rm a.s.}{\to}
 \frac{\xi^2_T}{1-2\alpha}$ (resp. $\xi^2_t\overset{\rm a.s.}{\to}\xi^2_T$)
as $t\to T$
by Lemma \ref{lm7} (resp. Lemma \ref{lm3}); hence, as desired one gets that $\alpha - \widehat{\alpha}_t\overset{\rm a.s.}{\to} 
0$ as $t\to T$.

When $\alpha=\frac12$, the identity (\ref{sothat}) becomes 
\begin{equation}\label{sothat2}
\alpha - \widehat{\alpha}_t =
\frac{\xi_t^2}{2\int_0^t \xi_u^2(T-u)^{-1}du
};
\end{equation}
as $t\to T$, we have $\int_0^t \xi_u^2(T-u)^{-1}du\overset{\rm a.s.}{\sim} |\log(T-t)|\xi^2_T$ 
(resp. $\xi^2_t\overset{\rm a.s.}{\to}\xi^2_T$) by Lemma \ref{lm7} (resp. Lemma \ref{lm3}).
Hence, here again  we have $\alpha - \widehat{\alpha}_t\overset{\rm a.s.}{\to}  0$ as $t\to T$.

Suppose now that $\alpha\in(\frac12,H)$.
As $t\to T$, we have $\int_0^t \xi_u^2(T-u)^{2\alpha-2}du \overset{\rm a.s.}{\to} \int_0^T \xi_u^2(T-u)^{2\alpha-2}du$ 
(resp. $\xi^2_t\overset{\rm a.s.}{\to}\xi^2_T$) by Lemma \ref{lm7} (resp. Lemma \ref{lm3}).
Hence (\ref{sothat}) yields this time that 
$\alpha - \widehat{\alpha}_t\overset{\rm a.s.}{\to} 
\alpha-\frac12$ as $t\to T$, that is
$\widehat{\alpha}_t\overset{\rm a.s.}{\to} 
\frac12$.

Assume finally that $\alpha\geq H$. By (\ref{isoint}), we have
\begin{eqnarray*}
E\big[(T-t)^{2\alpha-1}\xi^2_t\big]
&=&(T-t)^{2\alpha-1}\int_0^t du(T-u)^{-\alpha}\int_0^t dv(T-v)^{-\alpha}|v-u|^{2H-2}\\
&=&(T-t)^{2\alpha-1}\int_{T-t}^Tdu\,u^{-\alpha}\int_{T-t}^Tdv\,v^{-\alpha}|v-u|^{2H-2}\\
&=&(T-t)^{2H-1}\int_1^{\frac{T}{T-t}}du\,u^{-\alpha}\int_{1}^{\frac{T}{T-t}}dv\,v^{-\alpha}|v-u|^{2H-2}\\
&=&(T-t)^{2H-1}\int_{\frac{T-t}{T}}^1du\,u^{\alpha-2H}\int_{\frac{T-t}{T}}^1dv\,v^{\alpha-2H}|v-u|^{2H-2}\\
&\leq&(T-t)^{2H-1}\int_{\frac{T-t}{T}}^1du\,u^{2\alpha-2H-1}
\int_{0}^{\frac{1}u}dv\,v^{\alpha-2H}|v-1|^{2H-2}\\
&\leq&c(T-t)^{2H-1}\int_{\frac{T-t}{T}}^1du\,u^{2\alpha-2H-1}
\times \left\{
\begin{array}{ll}
1&\,\,\mbox{if $\alpha<1$}\\
\\
1+|\log u|&\,\,\mbox{if $\alpha=1$}\\
\\
u^{1-\alpha}&\,\,\mbox{if $\alpha>1$}
\end{array}
\right.
\\
&\leq&c(T-t)^{2H-1}\times \left\{
\begin{array}{ll}
|\log(T-t)|&\,\,\mbox{if $\alpha=H$}\\
\\
1&\,\,\mbox{if $\alpha>H$}
\end{array}
\right.\\
&\longrightarrow& 0\mbox{ as $t\to T$}.
\end{eqnarray*}
Hence, having a look at (\ref{sothat}) and because $\int_0^t \xi^2_u(T-u)^{2\alpha-2}du\overset{\rm a.s.}{\to} \int_0^T \xi^2_u(T-u)^{2\alpha-2}du
\in(0,\infty]$ as $t\to T$, we deduce that $\alpha - \widehat{\alpha}_t\overset{\rm prob.}{\to} 
\alpha-\frac12$ as $t\to T$, that is
$\widehat{\alpha}_t\overset{\rm prob.}{\to} 
\frac12$.

The proof of Theorem \ref{prop7} is done.
\qed

\bigskip
\noindent
{\it Proof of Theorem \ref{prop8}}. 1. Assume that $\alpha$ belongs to $(0,1-H)$.
We have, by using Lemma \ref{lm2} to go from the first to the second line,
\begin{eqnarray}
&&(T-t)^{\alpha-H}\big( \alpha - \widehat{\alpha}_t\big)  
=\frac{(T-t)^{1-H-\alpha}\eta_t}{(T-t)^{1-2\alpha}\int_0^t \xi_s^2(T-s)^{2\alpha-2}ds}\notag\\
&=&\frac{(T-t)^{1-H-\alpha}
\int_0^t (T-u)^{\alpha-1} dB_u\int_0^t (T-s)^{-\alpha} dB_s 
}
{(T-t)^{1-2\alpha}\int_0^t \xi_s^2(T-s)^{2\alpha-2}ds}
\notag\\
&&
-
\frac{(T-t)^{1-H-\alpha}\int_0^t \delta B_s(T-s)^{-\alpha}\int_0^s \delta B_u(T-u)^{\alpha-1}}
{(T-t)^{1-2\alpha}\int_0^t \xi_s^2(T-s)^{2\alpha-2}ds}\notag\\
&&-H(2H-1)\frac{(T-t)^{1-H-\alpha}\int_0^t ds\,(T-s)^{-\alpha}\int_0^s du\,(T-u)^{\alpha-1}(s-u)^{2H-2}}
{(T-t)^{1-2\alpha}\int_0^t \xi_s^2(T-s)^{2\alpha-2}ds}\notag
\end{eqnarray}
\begin{eqnarray}
&=&\frac{1-2\alpha}{\xi_T}\,
(T-t)^{1-H-\alpha} \int_0^t (T-u)^{\alpha-1} dB_u\notag\\
&&\hskip2cm
\times
\frac{ \int_0^t (T-s)^{-\alpha}dB_s  }{\xi_T}
\times
\frac{\xi_T^2}
{(1-2\alpha)(T-t)^{1-2\alpha}\int_0^t \xi_s^2(T-s)^{2\alpha-2}ds}
\notag\\
&&
-
\frac{(T-t)^{1-H-\alpha}\int_0^t \delta B_s(T-s)^{-\alpha}\int_0^s \delta B_u(T-u)^{\alpha-1}}
{(T-t)^{1-2\alpha}\int_0^t \xi_s^2(T-s)^{2\alpha-2}ds}\notag\\
&&-H(2H-1)\frac{(T-t)^{1-H-\alpha}\int_0^t ds\,(T-s)^{-\alpha}\int_0^s du\,(T-u)^{\alpha-1}(s-u)^{2H-2}}
{(T-t)^{1-2\alpha}\int_0^t \xi_s^2(T-s)^{2\alpha-2}ds}\notag\\
&=&a_t \times b_t\times c_t - d_t -e_t,\label{abcde}
\end{eqnarray}
with clear definitions for $a_t$ to $e_t$.
Lemma \ref{lm4} yields  
\[
a_t\overset{\rm law}{\to} (1-2\alpha)\sqrt{H(2H-1)\frac{\beta(2-\alpha-2H,2H-1)}{1-H-\alpha}}\times\frac{G}{\xi_T}\quad\mbox{as $t\to T$,}
\]
where $G\sim\mathcal{N}(0,1)$ is independent of $B$, whereas
Lemma \ref{lm3} (resp. Lemma \ref{lm7}) implies that $b_t\overset{\rm a.s.}{\to} 1$ (resp. $c_t\overset{\rm a.s.}{\to} 1$) as $t\to T$.
On the other hand, by combining Lemma \ref{lm7} with Lemma \ref{lm6} (resp. Lemma \ref{lm1}), we deduce that 
$d_t\overset{\rm prob.}{\to} 0$ (resp. $e_t\overset{\rm prob.}{\to} 0$) as $t\to T$.
By plugging all these convergences together we get that, as $t\to T$,
\[
(T-t)^{\alpha-H}\big(\widehat{\alpha}_t - \alpha\big)\overset{\rm law}{\to} 
(1-2\alpha)\sqrt{H(2H-1)\frac{\beta(2-\alpha-2H,2H-1)}{1-H-\alpha}}\times
\frac{G}{\xi_T}.
\]
Because it is well-known that the ratio of two independent $\mathcal{N}(0,1)$-random variables is $\mathcal{C}(1)$-distributed, 
to conclude it remains to compute the variance $\sigma^2$ of $\xi_T\sim\mathcal{N}(0,\sigma^2)$. By (\ref{isoint}), we have:
\begin{eqnarray}
E[\xi^2_T]&=&H(2H-1)\int_0^T du(T-u)^{-\alpha}\int_0^Tdv(T-v)^{-\alpha}|v-u|^{2H-2}\notag\\
&=&H(2H-1)\int_0^T du\,u^{-\alpha}\int_0^Tdv\,v^{-\alpha}|v-u|^{2H-2}\notag\\
&=&2H(2H-1)\int_0^T du\,u^{-\alpha}\int_0^u dv\,v^{-\alpha}(u-v)^{2H-2}\notag\\
&=&2H(2H-1)\int_0^T u^{2H-2\alpha-1}du\int_0^1 v^{-\alpha}(1-v)^{2H-2}dv\notag\\
&=&
\frac{H(2H-1)}{H-\alpha}T^{2H-2\alpha}\beta(1-\alpha,2H-1),\label{rr}
\end{eqnarray}
and the proof of the first part of Theorem \ref{prop8} is done.\\

2. Assume that $\alpha=1-H$. The proof follows the same lines as in point 1 above. The counterpart of decomposition (\ref{abcde})
is here:
\begin{eqnarray}
\frac{(T-t)^{1- 2H}}{\sqrt{|\log(T-t)|}}\big( \alpha - \widehat{\alpha}_t\big) 
&=&\frac{2H-1}{\xi_T
\sqrt{|\log(T-t)|}
} \int_0^t (T-s)^{-H} dB_s\notag\\
&&
\times
\frac{ \int_0^t (T-u)^{H-1}dB_u  }{\xi_T}
\times
\frac{\xi_T^2}
{(2H-1)(T-t)^{2H-1}\int_0^t \xi_s^2(T-s)^{-2H}ds}
\notag
\end{eqnarray}
\begin{eqnarray}
&&
-
\frac{\int_0^t \delta B_s(T-s)^{H-1}\int_0^s \delta B_u(T-u)^{-H}}
{\sqrt{|\log(T-t)|}(T-t)^{2H-1}\int_0^t \xi_s^2(T-s)^{-2H}ds}\notag\\
&&-H(2H-1)\frac{\int_0^t ds\,(T-s)^{H-1}\int_0^s du\,(T-u)^{-H}(s-u)^{2H-2}}
{\sqrt{|\log(T-t)|}(T-t)^{2H-1}\int_0^t \xi_s^2(T-s)^{-2H}ds}\notag\\
&=&\widetilde{a}_t \times \widetilde{b}_t\times \widetilde{c}_t - \widetilde{d}_t -\widetilde{e}_t.\notag
\end{eqnarray}
Lemma \ref{lm4} yields  
\[
\widetilde{a}_t\overset{\rm law}{\to}
(2H-1)^{\frac32}\sqrt{2H\,\beta(1-H,2H-1)}\times\frac{G}{\xi_T}
\quad\mbox{as $t\to T$,}
\]
where $G\sim\mathcal{N}(0,1)$ is independent of $B$, whereas
Lemma \ref{lm3} (resp. Lemma \ref{lm7}) implies that $\widetilde{b}_t\overset{\rm a.s.}{\to} 1$ (resp. 
$\widetilde{c}_t\overset{\rm a.s.}{\to} 1$) as $t\to T$.
On the other hand, by combining Lemma \ref{lm7} with Lemma \ref{lm6} (resp. Lemma \ref{lm1}), we deduce that 
$\widetilde{d}_t\overset{\rm prob.}{\to} 0$ (resp. $\widetilde{e}_t\overset{\rm prob.}{\to} 0$) as $t\to T$.
By plugging all these convergences together we get that, as $t\to T$,
\[
\frac{(T-t)^{1- 2H}}{\sqrt{|\log(T-t)|}}
\big(\widehat{\alpha}_t - \alpha\big)
\overset{\rm law}{\to} 
(2H-1)^{\frac32}\sqrt{2H\,\beta(1-H,2H-1)}\times\frac{G}{\xi_T}.
\]
Moreover, by (\ref{rr}) we have that $\xi_T\sim\mathcal{N}\big(0,HT^{4H-2}\beta(H,2H-1)
\big)$.
Thus, 
\[
(2H-1)^{\frac32}\sqrt{2H\,\beta(1-H,2H-1)}\times\frac{G}{\xi_T}
\overset{\rm law}{=}T^{1-2H}
(2H-1)^{\frac32}\sqrt{\frac{2\,\beta(1-H,2H-1)}{\beta(H,2H-1)}}
\times\mathcal{C}(1),
\]
and the convergence in point 2 is shown.
\\

3. Assume that $\alpha$ belongs to $(1-H,\frac12)$. Using the decomposition
\begin{eqnarray*}
(T-t)^{2\alpha-1}\big( \alpha - \widehat{\alpha}_t\big)
=\frac{\eta_t}
{(T-t)^{1-2\alpha}\int_0^t \xi_u^2(T-u)^{2\alpha-2}du},
\end{eqnarray*}
we immediately see that the second part of Theorem \ref{prop8} is an obvious consequence of Lemmas \ref{lm3bis} and \ref{lm7}.\\

4. Assume that $\alpha=\frac12$. Recall the identity (\ref{sothat2}) for this particular value of $\alpha$: 
\[
\alpha - \widehat{\alpha}_t =
\frac{\xi_t^2}{2\int_0^t \xi_u^2(T-u)^{-1}du
}.
\]
As $t\to T$, we have $\xi^2_t\overset{\rm a.s.}{\to}\xi^2_T$ by Lemma \ref{lm3}, whereas
$\int_0^t \xi_u^2(T-u)^{-1}du\overset{\rm a.s.}{\sim} |\log(T-t)|\xi^2_T$ by Lemma \ref{lm7}.
Therefore,  we deduce as announced that $ |\log(T-t)|\big(\alpha - \widehat{\alpha}_t\big)\overset{\rm a.s.}{\to}  \frac12$ as $t\to T$.
\qed

\vskip1cm

\noindent
{\bf Acknowledgments}. We thank an anonymous referee for his/her careful reading of the manuscript and for his/her valuable suggestions
and remarks.
We also thank Jingqi Han for pointing out a misprint in one of the estimates of the proof  
of Lemma \ref{lm3bis}.

\bibliographystyle{amsplain}

\end{document}